\documentclass[11pt]{article}

\usepackage{graphicx}
\usepackage{amsmath}
\usepackage{amsfonts}
\usepackage{fancyhdr} 
\usepackage{mathrsfs}
\usepackage{wrapfig}
\usepackage{bbold}

\newtheorem{thm}{Theorem}
\newtheorem{lem}{Lemma}

\newtheorem{prop}[lem]{Proposition}

\newcommand{\eee}{\mathcal{E}}

\newcommand{\sss}{\mathcal{S}}

\newcommand{\real}{\mathbb{R}}
\newcommand{\complex}{\mathbb{C}}

\newcommand{\ddd}{\mathbb{D}}

\newcommand{\so}{\mathcal{S}(\Omega)}
\newcommand{\sd}{\mathcal{S}(\ddd)}
\newcommand{\sso}{\mathscr{S}(\Omega)}

\newcommand{\an}{\mathfrak{A}}
\newcommand{\sao}{\mathcal{S}_\mathfrak{A}(\Omega)}
\newcommand{\ssao}{\mathscr{S}_\mathfrak{A}(\Omega)}
\newcommand{\ooo}{\mathcal{O}}

\newenvironment{pf}{\noindent {\em Proof}.\ \ }{\hspace*{\fill}\rule{.5ex}{1.4ex}\,}

\title{\vspace*{-1in}Constructive solutions to P\'olya-Schur problems}
\author{Peter C.~Gibson\footnote{Corresponding author} \footnote{Dept.~of Mathematics \& Statistics, York University, $\mathtt{pcgibson@yorku.ca}$} and Michael P.~Lamoureux\footnote{Dept.~of Mathematics \& Statistics, University of Calgary, $\mathtt{mikel@ucalgary.ca}$}}
\date{August 28, 2015}
\makeatletter
\let\newtitle\@title
\let\newauthor\@author
\let\newdate\@date
\makeatother

\begin{document}

\maketitle
\begin{abstract}
The general P\'olya-Schur problem is to characterize linear operators on the space of univariate polynomials that preserve stability, where a polynomial is stable with respect to a region $\Omega$ in the complex plane if it has no zeros in $\Omega$.  Stable preserving operators have proven to be important in a variety of applications ranging from statistical mechanics to combinatorics, and variants of P\'olya-Schur problems involving analytic functions are important in applications to signal processing.  We present a structure theorem that bridges polynomial and analytic P\'olya-Schur problems, providing constructive characterizations of stable-preserving operators for a general class of domains $\Omega$.  The structure theorem facilitates the solution of open P\'olya-Schur problems in the classical setting, and provides constructive characterizations of stable preserving operators in cases where previously known characterizations are non-constructive.  In the analytic setting, the structure theorem enables the explicit characterization of minimum-phase preserving operators on the half-line, a problem of importance in geophysical signal processing.   
\end{abstract}

\section{Introduction\label{sec-introduction}}

A univariate polynomial $p\in\complex[z]$ is said to be stable with respect to a set $\Omega\subset\complex$ if it has no zeros in $\Omega$. Let $\so$ denote the set of all $\Omega$-stable polynomials, together with 0.  A linear mapping
\[
A:\complex[z]\rightarrow\complex[z]
\]
is stability-preserving with respect to $\Omega$ if $A(\so)\subset\so$;  we denote the semigroup of all such mappings by $\sso$.   

The general P\'olya-Schur problem is to characterize $\sso$ for various $\Omega\subset\complex$; the name harkens back to work by P\'olya and Schur on multiplier sequences in the case $\Omega=\complex\setminus\real$, \cite{PoSc:1914}.  Recently, P\'olya-Schur problems have come into prominence for a variety of applications (see \cite{Wa:2011} for an overview), notably in the Ising model of statistical mechanics, where they underpin construction of Lee-Yang polynomials \cite{BoBr:2009Co,Ru:2010}.  These results rely on characterizations of $\sso$ for circular domains $\Omega$ and their boundaries presented in \cite{BoBr:2009An}.  The latter work gives two types of characterizations: algebraic, in which a particular class of test functions---such as translates of monomials---determines whether a given operator preserves stability; and transcendental, whereby an operator preserves stability according to whether its characteristic function belongs to a certain analytic class.  

Independently of the above developments, we have studied an analytic version of a P\'olya-Schur problem in the context of the Hardy-Hilbert space $H^2=H^2(\ddd)$ on the unit disk, motivated by geophysical applications \cite{GiLaMa:2011}.  More precisely, the classical factorization theorem for $H^2$ expresses an arbitrary function as a product of an inner function, defined to have constant modulus 1 on the boundary circle, and an outer function, by definition a cyclic vector of the unilateral shift (see \cite{MaRo:2007}).   Thus a function $f\in H^2$ is outer if and only if the span of functions of the form $z^nf(z)$, where $n\geq 0$, is dense in $H^2$.  Referring to functions of the form $z^nf(z)$ as shifted outer functions, \cite{GiLaMa:2011} constructively characterizes continuous linear operators $A:H^2\rightarrow H^2$ that preserve the class $\widetilde{\ooo}\subset H^2$ of all shifted outer functions.  This is analogous to a P\'olya-Schur problem because outer functions have no zeros in the open unit disk $\ddd$, and hence shifted outer functions have no zeros in the punctured disk $\ddd\setminus\{0\}$.  The characterization of operators that preserve shifted outer functions translates via the inverse $Z$-transform to a characterization of operators acting on causal digital signals that preserve the class of delayed minimum-phase signals, a property of importance for seismic recordings (see \cite{GiLaMa:2011} for full details).   

The present paper stems from two principal objectives: (i) to bring the classical P\'olya-Schur problems and their analytic analogues together within a common framework; and (ii) to extend the characterization of operators preserving delayed minimum-phase digital signals to operators acting on continuous signals in $L^2(\real_+)$, where the methods of \cite{GiLaMa:2011} are insufficient.  
We prove a structure theorem that achieves both of these objectives and that furthermore solves an array of previously open P\'olya-Schur problems.   
Let $\an(\ddd)$ denote the vector space of analytic functions on the open unit disk $\ddd$.   Given a set $\Omega\subset\ddd$, let $\sao\subset\an(\ddd)$ denote the family of functions having no zeros in $\Omega$, together with $0$. Let $\ssao$ denote the set of linear operators 
\[
A:\complex[z]\rightarrow\an(\ddd)\quad\mbox{ such that }\quad A(\sd)\subset\sao.
\]
As with polynomials, an analytic function in $\sao$ is said to be stable with respect to $\Omega$. 
Our main result characterizes the set $\ssao$ in terms of multiplication and composition operators, defined by the respective formulas
\[
M_\psi p=\psi p\quad\mbox{ and }\quad C_\varphi p=p\circ\varphi,\quad\mbox{ where }p\in\complex[z]\mbox{ and }\psi,\varphi\in\mathfrak{A}(\ddd).
\]
\begin{thm}[Structure Theorem]\label{thm-Structure}
Given a non-empty connected open set  $\Omega\subset\ddd$,  a linear map $A:\complex[z]\rightarrow\an(\ddd)$ has the property that 
$
A(\sd)\subset\sao
$
if and only if either:
\begin{enumerate}
\item there exist a function $\psi\in\sao$ and a linear functional $\nu:\complex[z]\rightarrow\complex$ such that $A(f)=\nu(f)\psi$, for all $f\in\complex[z]$; or 
\item there exist a function $\psi\in\sao$ and a non-constant function $\varphi\in\an(\ddd)$, where $\varphi(\Omega)\subset\ddd$, such that $A=M_\psi C_\varphi$.  
\end{enumerate}
\end{thm}

\pagestyle{fancyplain}

This result provides a bridge between the polynomial and analytic contexts. It applies to operators preserving any class $X$ of functions that are stable with respect to an appropriate region $\Omega$ provided only that 
\[
\so\subset X\subset\sao.
\] 
In \S\ref{sec-general} we use Theorem~\ref{thm-Structure} to solve the general P\'olya-Schur problem of characterizing linear operators $A:\complex[z]\rightarrow\complex[z]$ such that $A(\sss(\Omega_1))\subset\sss(\Omega_2)$, where $\Omega_1$ is bounded and $\Omega_2$ has non-empty interior (Theorem~\ref{thm-General}).   This solves one of the open problems posed by Borcea and Br\"and\'en in \cite{BoBr:2009An}; just as importantly, it gives a constructive solution in cases where there was a known, non-constructive characterization of $\sso$.  Furthermore, it removes altogether the previous restriction to circular domains $\Omega$.  

In \S\ref{sec-outer} we apply Theorem~\ref{thm-Structure} to the context of outer functions on the Hardy-Hilbert space $H^2$.  We characterize operators $A:H^2\rightarrow H^2$ that preserve the class of outer functions (Theorem~\ref{thm-Outer-preserving}), as well as those that preserve the class of shifted outer functions (Theorem~\ref{thm-Shifted-outer-preserving}).  While the latter result was proved earlier in \cite{GiLaMa:2011}, the method of proof used there depends essentially on the fact that moment functions $z^n$ are shifted outer functions.  We have included an alternate proof using Theorem~\ref{thm-Structure}, with $\Omega=\ddd\setminus\{0\}$, to show that the latter result subsumes the earlier methods.   It turns out that in analyzing operators on the space $L^2(\real_+)$ of causal signals with respect to a continuous time variable, the class of delayed minimum phase signals corresponds to modified outer functions with respect to Hardy space on the half-plane (see \cite{Ho:1962}).  This in turn can be related to $H^2$, but it is no longer shifted outer functions that come into play. 
Full details of the connection between minimum-phase-preserving operators on $L^2(\real_+)$ and outer-preserving operators on $H^2$ are given in \cite{GiLa:IP2012}, which is based on an earlier preprint version of the present paper.   

Our main result, Theorem~\ref{thm-Structure}, is proved in \S\ref{sec-proof}; the paper concludes with a brief discussion in \S\ref{sec-conclusions}.   

\section{Applications of the structure theorem\label{sec-applications}}

In this section we apply Theorem~\ref{thm-Structure} to both the polynomial and analytic contexts, obtaining new results in both.   In particular we solve the open problem (a) posed in \cite[\S4]{BoBr:2009An}, which corresponds to $\Omega_1=\Omega_2=\overline{\ddd}$ in Theorem~\ref{thm-General} below.  (Note that our notation differs from that of \cite{BoBr:2009An} in that they refer to the complement $\complex\setminus\Omega$, as opposed to the zero-free region $\Omega$.)   More than this, we give a constructive solution that is different from---and complementary to---existing algebraic or transcendental characterizations.

\subsection{General P\'olya-Schur problems\label{sec-general}}

\begin{thm}\label{thm-General}
Suppose $\Omega_1\subset\complex$ is bounded, and $\Omega_2\subset\complex$ has non-empty interior.  A linear map
$
A:\complex[z]\rightarrow\complex[z]
$
has the property that 
$
A(\sss(\Omega_1))\subset\sss(\Omega_2)
$
if and only if either:
\begin{enumerate}
\item there exist a linear functional 
$\nu:\complex[z]\rightarrow\complex$ and a polynomial $\psi\in\sss(\Omega_2)$ such that $A(f)=\nu(f)\psi$, for all $f\in\complex[z]$; or
\item there exist $\psi\in\sss(\Omega_2)$ and a non-constant polynomial $\varphi$ for which $\varphi(\Omega_2)\subset\Omega_1$ such that $A=M_\psi C_\varphi$.  
\end{enumerate}
\end{thm}

\begin{pf}
Given $\tau>0$, let $D_\tau:\complex[z]\rightarrow\complex[z]$ denote the dilation operator defined by the formula
\[
(D_\tau p)(z)=p(\tau z)\mbox{ for every } p\in\complex[z]\mbox{ and every }z\in\complex.
\]
Note that for any $\Omega\subset\complex$, 
$
D_\tau(\sss(\Omega))=\sss(\frac{1}{\tau}\Omega).
$

If $A:\complex[z]\rightarrow\complex[z]$ is a linear operator  such that 
$
A(\sss(\Omega_1))\subset\sss(\Omega_2),
$
where $\Omega_1$ is bounded and $\Omega_2$ has non-empty interior, then there exist $\delta>0$ such that $\delta\Omega_1\subset\ddd$, and $\varepsilon>0$ and an open connected set $\Omega\subset\complex$ such that 
\[
\Omega\subset\ddd\cap\varepsilon\Omega_2.
\]
It follows that for $\widetilde{A}=D_{1/\varepsilon}AD_\delta$, 
\[
\widetilde{A}\bigl(\sss(\ddd)\bigr)\subset\sss(\Omega).
\]
Thus by Theorem~\ref{thm-Structure}, either: 
\begin{enumerate}
\item there exist a function $\widetilde{\psi}\in\sao$ and a linear functional $\widetilde{\nu}:\complex[z]\rightarrow\complex$ such that $\widetilde{A}(f)=\widetilde{\nu}(f)\widetilde{\psi}$, for all $f\in\complex[z]$; or 
\item there exist a function $\widetilde{\psi}\in\sao$ and a non-constant function $\widetilde{\varphi}\in\mathfrak{A}(\ddd)$, where $\widetilde{\varphi}(\Omega)\subset\ddd$, such that $\widetilde{A}=M_{\widetilde{\psi}} C_{\widetilde{\varphi}}$.  
\end{enumerate}
In the former case, 
\[
A(f)=\left(\widetilde{\nu}D_{1/\delta}(f)\right)D_\varepsilon \widetilde{\psi}=\nu(f)\psi, 
\]
where $\psi=D_\varepsilon\widetilde{\psi}$ and $\nu=\widetilde{\nu}D_{1/\delta}$.  And in the latter case,
\[
A=D_\varepsilon M_{\widetilde{\psi}}C_{\widetilde{\varphi}}D_{1/\delta}=M_\psi C_\varphi,
\]
where $\psi=D_\varepsilon\widetilde{\psi}$ and $\varphi=\frac{1}{\delta}D_\varepsilon\widetilde{\varphi}$.  

It remains to verify that $\psi\in\sss(\Omega_2)$ and that the non-constant function $\varphi$ is a polynomial such that $\varphi(\Omega_2)\subset\Omega_1$; this is rather straightforward.  If $A$ is the zero operator then there is nothing to prove, since $\nu$ can be taken to be the zero functional, and $\psi$ can be taken arbitrarily; in particular one may take $\psi\in\sss(\Omega_2)$.  

If $A$ is non-zero of rank 1, then, by hypothesis, $A1=\nu(1)\psi\in\sss(\Omega_2)$, since $1\in\sss(\Omega_1)$.  Either $\nu(1)\neq0$---in which case $\nu(1)\psi\in\sss(\Omega_2)$ implies that $\psi\in\sss(\Omega_2)$---or $\nu(1)=0$ and there exists $n>0$ such that $\nu(z^n)\neq 0$.  In the latter case, fix $r>0$ sufficiently large that $\Omega_1\subset r\ddd$.  Then $q(z)=z^n-r^n\in\sss(\Omega_1)$, whereby $Aq=\nu(q)\psi\in\sss(\Omega_2)$, and $\nu(q)=\nu(z^n)\neq0$.  It follows that $\psi\in\sss(\Omega_2)$.  Thus $\psi\in\sss(\Omega_2)$ whatever the value of $\nu(1)$.

Otherwise $A=M_\psi C_\varphi$ has rank at least 2, and $A1=\psi\in\sss(\Omega_2)$.   To see that $\varphi$ is a polynomial such that
$
\varphi(\Omega_2)\subset\Omega_1,
$
let $z_0\in\complex\setminus\Omega_1$, so that $p(z)=z-z_0$ belongs to $\sss(\Omega_1)$.  Then $Ap\in\sss(\Omega_2)$, where
\[
(Ap)(z)=\psi(z)\varphi(z)-z_0\psi(z).
\]
Note that because $A$ has rank at least 2, $\varphi$ is not constant, and so $Ap$ cannot be 0.   Therefore,
\[
Ap=\psi(\varphi-z_0)\in\sss(\Omega_2)\setminus\{0\},
\]
which---given that $\psi\in\sss(\Omega_2)$---implies that for all $z\in\Omega_2$, $\varphi(z)\neq z_0$.   Since $z_0\in\complex\setminus\Omega_1$ was arbitrary, it follows that for every $z\in\Omega_2$, $\varphi(z)\in\Omega_1$, as desired.  That $\varphi\in\mathfrak{A}(\ddd)$ is actually a polynomial follows from the fact that $Az^n=\psi\varphi^n$ is a polynomial for each $n\geq0$.   

The converse direction of the theorem is much simpler.   If $A=M_\psi C_\varphi$ where $\psi\in\sss(\Omega_2)$ and $\varphi$ is a polynomial such that $\varphi(\Omega_2)\subset\Omega_1$, then it follows immediately that 
\[
A(\sss(\Omega_1))\subset\sss(\Omega_2),
\]
and similarly if $A(f)=\nu(f)\psi$ for every $f\in\complex[z]$, where $\nu:\complex[z]\rightarrow\complex$ is a linear functional.  \end{pf}

\subsection{Outer preserving operators on Hardy space\label{sec-outer}}

Next we consider the analytic context of the Hardy-Hilbert space $H^2=H^2(\ddd)$.   See \cite{GiLaMa:2011,GiLa:IP2012} for details on how this relates to signal processing and, in particular, to geophysics.
 
Let $\ooo$ denote the set of all outer functions in $H^2$, and let $\widetilde{\ooo}$ denote the set of all shifted outer functions.  
\begin{lem}\label{lem-Rank-1}
A bounded linear functional $\nu:H^2\rightarrow\complex$ satisfies
\[
\nu(\ooo)\subset\complex\setminus\{0\} \rule{10pt}{0pt} (\mbox{respectively, }\nu(\widetilde{\ooo})\subset\complex\setminus\{0\})
\]
if and only if there exist a point $z_0\in\ddd$ (respectively, $z_0\in\ddd\setminus\{0\}$) and a scalar $\sigma\in\complex\setminus\{0\}$ such that for all $f\in H^2$,
\[
\nu(f)=\sigma f(z_0).
\]
\end{lem}

\begin{pf}
If there exist a point $z_0\in\ddd$ and a scalar $\sigma\in\complex\setminus\{0\}$ such that for all $f\in H^2$,
\[
\nu(f)=\sigma f(z_0),
\]
then $\nu:H^2\rightarrow\complex$ is bounded (as evaluation at a point in $\ddd$ is a bounded linear functional on $H^2$) and if $f\in\ooo$ then $\sigma f(z_0)\neq0$, since outer functions have no zeros in the interior of the unit disk.  Similarly for the case of shifted outer functions $\widetilde{\ooo}$ with $z_0\in\ddd\setminus\{0\}$.   

For the converse implication, suppose that $\nu:H^2\rightarrow\complex$ is a bounded linear functional such that 
\[
\nu(\ooo)\subset\complex\setminus\{0\}.
\]
Set $\rho_n=\nu(z^n)$ for each $n\geq0$; and for each $\xi\in\ddd$, let $g_\xi\in H^2$ denote the function
\[
g_\xi(z)=\sum_{n=0}^\infty (\xi z)^n.
\]
Since $\nu$ is well-defined and bounded, the series 
\[
\nu (g_\xi)=\sum_{n=0}^\infty \rho_n\xi^n
\]
converges for every $\xi\in\ddd$.  In particular, the inequality $|\rho_n|<r^{-n}$ holds eventually for each fixed $0<r<1$.  For each $w\in\complex$, let $f_w$ denote the scaled exponential 
\[
f_w(z)=e^{wz}.
\]
Each $f_w\in\ooo$, and therefore by hypothesis the function $F(w)$ defined as
\[
F(w)=\nu(f_w)=\sum_{n=0}^\infty \rho_n\frac{w^n}{n!}
\]
is zero-free.  Moreover, the eventual inequality $|\rho_n|<r^{-n}$ shows $F$ to be entire of order at most 1, whence  
\[
F(w)=e^{\alpha+\beta w},
\]
for some scalars $\alpha$ and $\beta$, by Hadamard's Theorem.  It follows that 
\[
\rho_n=e^\alpha \beta^n,
\]
 for each $n\geq0$.  Furthermore, the Riesz representation theorem for linear functionals implies that there exists a $g\in H^2$ such that 
\[
\nu(f)=\langle f,g\rangle \mbox{ for every }f\in H^2.
\]
The coefficients of $g(z)=\sum_{n=0}^\infty a_nz^n$ are given by 
\[
\overline{a}_n=\nu(z^n)=\rho_n=e^\alpha\beta^n;
\]
therefore $g$ has the power series expansion
\[
g(z)=e^\alpha\sum_{n=0}^\infty\overline{\beta}^nz^n,
\]
and $|\beta|<1$ since $g\in H^2$. Thus $\langle f,g\rangle=e^\alpha f(\beta)$.  Setting $z_0=\beta$ and $\sigma=e^\alpha$, this shows that 
\[
\nu(f)=\sigma f(z_0),
\]
as desired.  The case where $\nu(\widetilde{\ooo})\subset\complex\setminus\{0\}$ is similar, except that the fact that $z\in\widetilde{\ooo}$ implies that $\nu(z)=e^\alpha\beta\neq0$, whereby $\beta\in\ddd\setminus\{0\}$.   
\end{pf}

\begin{thm}\label{thm-Outer-preserving}
Let $A:H^2\rightarrow H^2$ be a bounded linear operator such that 
\[
A(\ooo)\subset\ooo.
\]
Then there exist an analytic function $\varphi:\ddd\rightarrow\ddd$ and a function $\psi\in\ooo$ such that 
\[
A=M_\psi C_\varphi.
\] 
\end{thm}

\begin{pf}
Since $\complex[z]\subset H^2\subset\mathfrak{A}(\ddd)$ and $\sss(\ddd)\subset\ooo\subset\sss_{\mathfrak{A}}(\ddd)$, the hypothesis of the theorem implies that the restriction of $A$ to polynomials is a linear operator of the form 
\[
A:\complex[z]\rightarrow\mathfrak{A}(\ddd)
\]
such that $A(\sss(\ddd))\subset\sss_{\mathfrak{A}}(\ddd)$.  Thus Theorem~\ref{thm-Structure} applies with $\Omega=\ddd$, and hence either:
\begin{enumerate}
\item there exist a function $\widetilde{\psi}\in\sss(\ddd)$ and a linear functional $\nu:\complex[z]\rightarrow\complex$ such that $A(f)=\nu(f)\widetilde{\psi}$ for every $f\in\complex[z]$; or 
\item there exist a function $\psi\in\sss(\ddd)$ and a non-constant analytic function $\varphi:\ddd\rightarrow\ddd$, such that $A=M_\psi C_\varphi$.  
\end{enumerate}
In case~1, the hypothesis that $\nu(\ooo)\widetilde{\psi}\subset(\ooo)$ implies that $\widetilde{\psi}\in\ooo$ and $\nu(\ooo)\subset\complex\setminus\{0\}$.   Since $A$ is bounded, this implies in turn that $\nu:H^2\rightarrow\complex$ is a bounded linear functional, since $||\widetilde{\psi}\nu(f)||=||\widetilde{\psi}||\,|\nu(f)|$ for each $f\in H^2$.   Therefore, by Lemma~\ref{lem-Rank-1}, the linear functional $\nu$ is proportional to evaluation at a point: for every $f\in H^2$, 
\[
\nu(f)=\sigma f(z_0), \mbox{ for some }z_0\in\ddd \mbox{ and some non-zero }\sigma\in\complex.
\]
Thus, letting $\varphi:\ddd\rightarrow\ddd$ denote the constant function $\varphi(z)=z_0$, and setting $\psi=\sigma\widetilde{\psi}$, the rank~1 operator $A$ has the form $A=M_\psi C_\varphi$, in conformity with the conclusion of the present theorem.  

In case~2, where $A$ has the form $M_\psi C_\varphi$, the fact that $1\in\ooo$ implies that $\psi=A1\in\ooo$, giving the desired conclusion once again.  \end{pf}

The following fact about shifted outer functions, which is needed for the proof of Theorem~\ref{thm-Shifted-outer-preserving},  follows easily from the standard integral representation for outer functions (see \cite{MaRo:2007}).

\begin{prop}\label{prop-Shifted-product}
Let $f,g\in H^2$.  If $f\in\widetilde{\ooo}$ and $fg\in\widetilde{\ooo}$, then $g\in\widetilde{\ooo}$.
\end{prop}

\begin{thm}[See \cite{GiLaMa:2011}]\label{thm-Shifted-outer-preserving}
Let $A:H^2\rightarrow H^2$ be a bounded linear operator such that 
\[
A\left(\widetilde{\ooo}\right)\subset\widetilde{\ooo}.
\]
Then there exist functions $\psi,\varphi\in\widetilde{\ooo}$, where $\varphi:\ddd\rightarrow\ddd$, such that 
\[
A=M_\psi C_\varphi.
\] 
\end{thm}

\begin{pf}
The proof is very similar to that of Theorem~\ref{thm-Outer-preserving}.  Setting $\Omega=\ddd\setminus\{0\}$, note that $\sss(\ddd)\subset\widetilde{\ooo}\subset\sss_{\mathfrak{A}}(\Omega)$; thus the hypothesis of the theorem implies that 
\[
A:\complex[z]\rightarrow\mathfrak{A}(\ddd),
\]
with $A(\sss(\ddd))\subset\sss_{\mathfrak{A}}(\Omega)$.   Theorem~\ref{thm-Structure} therefore implies that either:
\begin{enumerate}
\item there exist a function $\widetilde{\psi}\in\sss(\ddd)$ and a linear functional $\nu:\complex[z]\rightarrow\complex$ such that $A(f)=\nu(f)\widetilde{\psi}$ for all $f\in\complex[z]$; or 
\item there exist a function $\psi\in\sss(\ddd)$ and a non-constant function $\varphi\in\mathfrak{A}(\ddd)$, where $\varphi(\Omega)\subset\ddd$, such that $A=M_\psi C_\varphi$.  
\end{enumerate}
In case~1, as in the proof of Theorem~\ref{thm-Outer-preserving} above,  Lemma~\ref{lem-Rank-1} yields that the linear functional $\nu$ is proportional to evaluation at a point: for every $f\in H^2$, 
\[
\nu(f)=\sigma f(\zeta), \mbox{ for some }\zeta\in\Omega\mbox{ and some non-zero }\sigma\in\complex,
\]
with the difference that $\zeta\neq0$, as per the part of Lemma~\ref{lem-Rank-1} pertaining to $\widetilde{\ooo}$.  
As before, letting $\varphi:\ddd\rightarrow\ddd$ denote the constant function $\varphi(z)=\zeta$, and setting $\psi=\sigma\widetilde{\psi}$, the rank~1 operator $A$ has the form $A=M_\psi C_\varphi$, where the non-zero constant function $\varphi$ belongs to $\widetilde{\ooo}$.  

In case~2, where $A$ has the form $M_\psi C_\varphi$, the fact that $1\in\widetilde{\ooo}$ implies that $\psi=A1\in\widetilde{\ooo}$.   And the identity function $f(z)=z$ belongs to $\widetilde{\ooo}$, so $\psi\varphi=Af\in\widetilde{\ooo}$ also.  By Proposition~\ref{prop-Shifted-product} this implies that $\varphi\in\widetilde{\ooo}$.  Since $\varphi\in\mathfrak{A}(\ddd)$ is analytic, the property $\varphi(\Omega)\subset\ddd$ given by Theorem~\ref{thm-Structure} implies further that $\varphi(\ddd)\subset\ddd$, completing the proof. \end{pf}

\section{Proof of the structure theorem\label{sec-proof}}

Hadamard's theorem, \cite[Thm~8]{Al:1979} or \cite[Ch~I,\S10]{Le:1964}, applied to entire functions of exponential type provides a key analytic tool used to prove Theorem~\ref{thm-Structure}.  For present purposes an entire function $f:\complex\rightarrow\complex$ is of exponential type if $f(z)\leq Ce^{|z|}$ for some $C>0$; \cite{Le:1964} serves as an essential reference on the subject.  (Such functions are sometimes also referred to as having order 1; see \cite{Hi:1962}.) Hadamard's Theorem constrains the Weierstrass product form of what we call the characteristic functions (defined later in this section) of a given operator $A:\complex[z]\rightarrow\an(\ddd)$ preserving stability in the sense of Theorem~\ref{thm-Structure}.  This in turn allows us to determine the precise structure of the operator itself.  Our argument makes essential use of both the first and second characteristic functions.  Whereas (versions of) the first characteristic function appear elsewhere in the literature, it appears that the second characteristic function has not been previously recognized as important.  Before Hadamard's theorem can be applied, we need first to infer some basic consequences of the hypothesis that $A(\sd)\subset\sao$ on the moments of $A$.  

A linear map $A:\complex[z]\rightarrow\an(\ddd)$ is evidently determined by its action on monomials $z^n$; thus for each $n\geq0$, denote the $n$th moment of $A$ by
\[
\psi_n=A(z^n).
\]
For the remainder of the present section $\Omega$ denotes a connected open subset of $\ddd$, unless the contrary is explicitly stated.   

It turns out that for $A\in\ssao$, either all the moments $\psi_n$ lie on a single complex line, or else for each $z\in\Omega$, the numbers $|\psi_n(z)|$ are bounded uniformly in $n$.   This dichotomy is expressed in detail by the following two propositions. 
\begin{prop}\label{prop-zero-zeroeth-moment}
Suppose $A\in\ssao$ and that $\psi_0$ is identically zero.  Then there exists $\varphi\in\an(\ddd)$ such that $\psi_n\in\complex\varphi$ for every $n\geq0$; i.e., the operator $A$ has rank at most 1.  
\end{prop}
\begin{pf}
If $A$ is not the zero operator then for some $n\geq1$, the $n$th moment $\varphi=\psi_n$ is not identically zero. 
Since $\varphi$ is analytic and $\Omega$ contains a condensation point, there exists a point $\zeta\in\Omega$ at which $\varphi(\zeta)\neq0$.   Given an arbitrary $\psi_m$, there is a choice of $a\in\complex$ such that 
\[
\psi_m(\zeta)-a\varphi(\zeta)=0,
\]
whereby the function $\psi_m-a\varphi$ is not $\Omega$-stable.  Note that for any $\alpha\in\complex$ with $|\alpha|\geq|a|+1$, the polynomial 
\[
p(z)=\alpha+z^m-az^n
\]
is $\ddd$-stable.  Therefore since $A\in\ssao$, $Ap=\psi_m-a\varphi$ is either $\Omega$-stable or identically 0.  The former possibility has been ruled out, forcing $\psi_m=a\varphi$.  Thus each of the moments of $A$ lies in the line $\complex\varphi$. 
\end{pf}

If all its moments are scalar multiples of a fixed function $\varphi$, then the operator $A$ has rank 1 and has the form 
\[
Ap=\nu(p)\varphi,
\]
where $\nu:\complex[z]\rightarrow\complex$ is a linear functional.  This represents the degenerate case of the main structure theorem.  The following proposition addresses the non-degenerate case.    
\begin{prop}\label{prop-non-zero-zeroeth-moment}
Suppose $A\in\ssao$ and that $A$ has rank at least 2.  Then for each $n\geq1$, the $n$th moment of $A$ is subject to the following bounds. 
\begin{enumerate}
\item If $\psi_n\not\in\complex\psi_0$ then $\left|\displaystyle\frac{\psi_n(z)}{\psi_0(z)}\right|<1$ for every $z\in\Omega$.  
\item If $\psi_n\in\complex\psi_0$ then $\left|\displaystyle\frac{\psi_n}{\psi_0}\right|<3$. 
\end{enumerate}
\end{prop}
\begin{pf}
Since $A$ has rank greater than 1, it follows from Proposition~\ref{prop-zero-zeroeth-moment} that $\psi_0$ is not identically 0.  And $\psi_0$ is therefore $\Omega$-stable since $1\in\sd$ and $A\in\ssao$. 

Suppose that $\psi_n\not\in\complex\psi_0$, for some fixed $n\geq1$, and let $\zeta\in\Omega$ be arbitrary.  
Note that for $a\in\complex\setminus\ddd$ the polynomial $p(z)=a+z^n$ is $\ddd$-stable.  Therefore $Ap=a\psi_0+\psi_n$ is either $\Omega$-stable or identically zero. But $\psi_n\not\in\complex\psi_0$, so the latter possibility is ruled out.  Since the value $a=-\psi_n(\zeta)/\psi_0(\zeta)$ renders $Ap$ unstable, it follows that 
$-\psi_n(\zeta)/\psi_0(\zeta)\in\ddd$ and hence that $|\psi_n(\zeta)|<|\psi_0(\zeta)|$, proving part 1.  

Next suppose that $\psi_n=\alpha\psi_0$ for some $\alpha\in\complex$.  By Proposition~\ref{prop-zero-zeroeth-moment}, there exists an $m\geq1$ for which $\psi_m\not\in\complex\psi_0$, since $A$ has rank at least two.  Let $\zeta\in\Omega$, and set $a=-\psi_m(\zeta)/\psi_0(\zeta)$.  Then $|a|<1$, by part~1 above, and so the polynomial
\[
p(z)=(a-\alpha)+z^m+z^n
\]
is $\ddd$-stable as long as $|\alpha|\geq 3$.  But 
\[
Ap=a\psi_0-\alpha\psi_0+\psi_m+\psi_n=a\psi_0+\psi_m,
\]
which has a zero at $\zeta\in\Omega$ and is not identically zero.  Therefore $p$ cannot be stable, forcing $|\alpha|<3$. This proves part 2.\end{pf}

For integers $j\geq 1$, define the $j$th characteristic function of $A\in\ssao$ by the formula
\begin{equation}\label{char}
F_j(z,w)=\sum_{n=0}^\infty\psi_{nj}(z)\frac{w^n}{n!}.
\end{equation}
Note that the subscript $nj$ on the right-hand side of (\ref{char}) denotes a product (not a pair of indices); in the present paper only the cases $j=1,2$ are needed.   
\begin{prop}\label{prop-Fkentire}
The following statements hold for any $A\in\ssao$ of rank at least 2 and any $j\geq1$.
\begin{enumerate}
\item With respect to $w$, the characteristic function $F_j(z,w)$ is an entire function of exponential type; and $F_j(z,w)$ is analytic with respect to $z\in\Omega$. 
\item If $F_j(z_0,w_0)=0$ for some $(z_0,w_0)\in\Omega\times\complex$, then $F_j(z,w_0)=0$ for all $z\in\Omega$; moreover, the order of $w_0$ as a zero of $F(z,w)$ is independent of $z$.  
\end{enumerate} 
\end{prop}
\begin{pf}
For fixed $z \in \Omega$, the function $F_j(z,w)$ is entire in $w$ by the bound 
\[
|\psi_{jn}(z)|\leq 3|\psi_0(z)|.
\]
That it has order at most 1 follows from the estimate
\[
|F_j(z,w)|\leq 3|\psi_0(z)|e^{|w|}.
\]
For each compact subset $K\subset\Omega$, the function $\psi_0$ is bounded on $K$ whence the functions $\psi_n$ are uniformly bounded on $K$ by Proposition~\ref{prop-non-zero-zeroeth-moment}.  So for fixed $w$, the series (\ref{char}) is uniformly Cauchy on $K$, and thus the series converges to an analytic function in $z$, proving part 1. 

Fix $w_0\in\complex$ and set 
\[
\sigma_n(z)=\sum_{m=0}^n\frac{(w_0z^k)^m}{m!},
\]
the $n$th partial sum of the Taylor series expansion of $e^{w_0z^k}$.  
As $e^{w_0 z^k}$ is bounded away from zero on the closed disk, and the partial sums converge uniformly on this disk, the $\sigma_n$ are eventually $\ddd$-stable. The images $A(\sigma_n)$ are either $\Omega$-stable, or identically zero, and converge uniformly to $F_j(z, w_0)$.  If $F_j(z,w_0)$ is not identically zero, it can be written as a limit of a subsequence of $\Omega$-stable functions $ A(\sigma_n)$, and thus by Hurwitz's theorem, $F_j(z,w_0)$ has no zeros at all in $\Omega$.

To see that the order of each zero is independent of $z$, fix a point $w_0$ in the set 
\[
E_j=\left\{w\in\complex\,|\,F_j(z,w)=0\right\}.
\]
Choose a circular path $\gamma$ in $\complex$ centred at $w_0$ of sufficiently small radius that every other member of $E_j$ lies strictly outside $\gamma$.   Since $F_j(z,w)$ is analytic in $z$, the order of the zero $w_0$, given by the integer valued function 
\[
\xi(z)=\frac{1}{2\pi i}\int_{\gamma}\frac{\frac{\partial F_j}{\partial w}(z,w)}{F_j(z,w)}\,dw,
\]
is continuous, and therefore constant, since $\Omega$ is connected. 
\end{pf}

Based on Proposition~\ref{prop-Fkentire}, let $\eee_j$ denote the sequence of zeros distinct from 0
\[
\left\{w\in\complex\,|\,w\neq0\;\&\;F_j(z,w)=0\right\}
\]
listed according to multiplicity; thus $\eee_j$ may be empty, finite, or countably infinite.  Let $\nu_j\geq0$ denote the order of $w=0$ as a zero of  $F_j(z,w)$.  In what follows we adopt the convention that products indexed by the empty set denote the constant function 1.

\begin{prop}\label{prop-10mod}
Suppose that $A\in\ssao$ has rank at least 2.  Then the first two characteristic functions have the form
\[
\begin{split}
F_1(z,w)&=q_1(w)e^{\alpha(z)+w\beta_1(z)}\\
F_2(z,w)&=q_2(w)e^{\alpha(z)+w\beta_2(z)}
\end{split}
\]
where $\beta_1$ is non-constant, and where $q_1$ and $q_2$ are entire of the form
\[
\begin{split}
q_1(w)&=\sum_{n=0}^\infty c_nw^n\quad\mbox{ with }c_0=1\mbox{ and }c_1=0\\
q_2(w)&=\sum_{n=0}^\infty d_nw^n\quad\mbox{ with }d_0=1.
\end{split}
\]
\end{prop}

\begin{pf}  By Proposition~\ref{prop-Fkentire}, $F_j(z,w)$ is entire of order at most 1 as a function in $w$, for $j=1,2$.  Hadamard's Theorem implies that the genus of $F_j(z,w)$ as a function in $w$ is either 0 or 1.  If $F_j$ has genus 0 then 
\begin{equation}\label{reciprocalrootsum}
\sum_{w_n\in\eee_j}\frac{1}{|w_n|}<\infty
\end{equation}
and the canonical Weierstrass product representation of $F_j$ has the form 
\begin{equation}\label{genus0}
F_j(z,w)=e^{\alpha_j(z)}\,\,w^{\nu_j}\!\!\prod_{w_n\in\eee_j}\left(1-\frac{w}{w_n}\right).
\end{equation}
If $F_j$ has genus 1 then either (\ref{reciprocalrootsum}) holds and $F_j$ has canonical Weierstrass product of the form
\begin{equation}\label{genus1dummy}
F_j(z,w)=e^{\alpha_j(z)+\beta_j(z)w}\,\,w^{\nu_j}\!\!\prod_{w_n\in\eee_j}\left(1-\frac{w}{w_n}\right),
\end{equation}
where $\beta_j(z)\neq0$, or 
\[
\sum_{w_n\in\eee_j}\frac{1}{|w_n|}=\infty\quad\mbox{ and }\quad\sum_{w_n\in\eee_j}\frac{1}{|w_n|^2}<\infty,
\]
and $F_j$ has canonical Weierstrass product of the form
\begin{equation}\label{genus1}
F_j(z,w)=e^{\alpha_j(z)+\beta_j(z)w}\,\,w^{\nu_j}\!\!\prod_{w_n\in\eee_j}e^{\frac{w}{w_n}}\left(1-\frac{w}{w_n}\right).  
\end{equation}
(See \cite[Ch.~5,\S2.3]{Al:1979}.) Since (\ref{reciprocalrootsum}) implies
\[
\quad\sum_{w_n\in\eee_j}\frac{1}{|w_n|^2}<\infty,
\]
and hence that the product 
\[
\prod_{w_n\in\eee_j}e^{\frac{w}{w_n}}\left(1-\frac{w}{w_n}\right)
\]
converges uniformly on compact sets, if $F_j$ has canonical product (\ref{genus1dummy}), it may be represented also by the (non-canonical) formula
\begin{equation}\label{genus1alt}
F_j(z,w)=e^{\alpha_j(z)+(\beta_j(z)-\gamma_j)w}\,\,w^{\nu_j}\!\!\prod_{w_n\in\eee_j}e^{\frac{w}{w_n}}\left(1-\frac{w}{w_n}\right),
\end{equation}
where $\gamma_j=\sum_{w_n\in\eee_j}\frac{1}{|w_n|}$.  After relabeling this conforms to (\ref{genus1}); it is the representation that we shall analyze below.  Thus, to summarize, Hadamard's theorem implies that the characteristic function $F_j$ may be represented by a Weierstrass product of the form (\ref{genus0}) or (\ref{genus1}).  
Proposition~\ref{prop-Fkentire} also guarantees that the zeros $w_n$ (including their multiplicities) and the indices $\nu_j$ do not depend on $z$.  (Of course the values $\alpha_j$ and $\beta_j$ appearing in the Weierstrass products may depend on $z$.)  

In fact the index $\nu_j=0$ for $j=1,2$.  To see this, note that, by definition, $F_j(z,0)=\psi_0$.  Since $A$ has rank at least 2, it follows from Proposition~\ref{prop-zero-zeroeth-moment} that $\psi_0$ is not identically 0. According to the formulas (\ref{genus0}) and (\ref{genus1}), $F_j(z,0)=q_j(0)e^{\alpha(z)}$ which is non-zero only if $q_j(0)\neq0$ whereby $\nu_j=0$.   The given Weierstrass product formulas then yield that $q_j(0)=1$,  which proves that $c_0=d_0=1$, and that $\psi_0=e^{\alpha_j}$, from which it follows that $\alpha_1=\alpha_2$.   We thus drop the subscripts, writing $\alpha=\alpha_1=\alpha_2$. 

Comparing the Taylor series expansion of (\ref{genus0}) with respect to $w$ to the definition (\ref{char}) of the characteristic function $F_1$ shows the moments of $A$ all to be multiples of $e^{\alpha(z)}$, whereby $A$ has rank 1; therefore $F_1(z,w)$ has the form (\ref{genus1}) since $A$ has rank at least 2.   If $\beta_1(z)$ is constant, than $A$ is again seen to have rank 1, so $\beta_1(z)$ is non-constant.  Since the coefficient of $w$ in the Taylor expansion of each elementary factor 
\[
e^{\frac{w}{w_n}}\left(1-\frac{w}{w_n}\right)
\]
is 0, the coefficient $c_1$ of $w$ is necessarily 0 in the expansion of $q_1(w)$.  Allowing $\beta_2=0$, the formulation $F_2(z,w)=q_2(w)e^{\alpha(w)+w\beta_2(z)}$ covers both possibilities (\ref{genus0}) and (\ref{genus1}) for $j=2$.   
Hadamard's Theorem itself guarantees that each $q_j(w)$ is entire.  
\end{pf}

It will be convenient to normalize the moments of $A$ by dividing by $\psi_0$.  To this end we define $\varphi_n=\psi_n/\psi_0$ for each $n\geq0$.  With this notation we have that 
\begin{equation}\label{charmod}
\frac{F_j(z,w)}{\psi_0(z)}=\sum_{n=0}^\infty\varphi_{nj}(z)\frac{w^n}{n!}
\end{equation}
and also that 
\begin{align}
\frac{F_1(z,w)}{\psi_0(z)}&=q_1(w)e^{w\beta_1(z)}\nonumber\\
&=\left(\sum_{n=0}^\infty c_nw^n\right)\sum_{n=0}^\infty\beta_1(z)^n\frac{w^n}{n!}\nonumber\\
&=\sum_{n=0}^\infty n!\left(\sum_{j=0}^nc_{n-j}\frac{\beta_1(z)^j}{j!}\right)\frac{w^n}{n!}\label{F1}\\
\frac{F_2(z,w)}{\psi_0(z)}&=q_2(w)e^{w\beta_2(z)}\nonumber\\
&=\sum_{n=0}^\infty n!\left(\sum_{j=0}^nd_{n-j}\frac{\beta_2(z)^j}{j!}\right)\frac{w^n}{n!}\label{F2}
\end{align}

We will show that the only way (\ref{F1}) and (\ref{F2}) can be consistent with (\ref{charmod}) is if $q_1=1$ is constant.  
\begin{prop}\label{prop-q1}
Suppose that $A\in\ssao$ has rank at least 2.  Then the first characteristic function of $A$ has the form
\[
F_1(z,w)=e^{\alpha(z)+w\beta(z)}.
\]
\end{prop}
\begin{pf}
To begin, consider $\varphi_2(z)$.  According to (\ref{charmod}) this is the coefficient of $w^2/2$ in (\ref{F1}) and the coefficient of $w$ in (\ref{F2}).   Thus we have 
\[
\varphi_2(z)=2c_2+\beta_1(z)^2=d_1+\beta_2(z).
\]
It is not possible for $\varphi_2(z)=d_1+\beta_2(z)$ to be independent of $z$, since this contradicts the assertion in Proposition~\ref{prop-10mod} that $\beta_1(z)$ is non-constant.  This shows that $\beta_2$ is also not constant.  It follows further that $d_1=0$, just as $c_1=0$ in the proof of Proposition~\ref{prop-10mod}, since the Weierstrass product form for $F_2(z,w)$ involves only even powers of $w$ if $\beta_2$ is non-constant.  Thus
\[
\beta_2(z)=2c_2+\beta_1(z)^2.
\]

The next step is to consider $\varphi_4(z)$ in light of the facts that $d_1=0$ and $\beta_2(z)=2c_2+\beta_1(z)^2$.   On one hand, $\varphi_4$ is the coefficient of $w^4/4!$ in (\ref{F1}), so
\begin{equation}\label{phi-4-1}
\varphi_4(z)=4!(c_4+c_3\beta_1(z)+c_2\beta_1(z)^2/2+\beta_1(z)^4/4!).
\end{equation}
On the other hand, $\varphi_4$ is the coefficient of $w^2/2$ in (\ref{F2}), whereby
\begin{align}
\varphi_4(z)&=2(d_2+\beta_2(z)^2/2)\nonumber\\
&=2(d_2+(2c_2+\beta_1(z)^2)^2/2)\nonumber\\
&=2d_2+4c_2^2+4c_2\beta_1(z)^2+\beta_1(z)^4.\label{phi-4-2}
\end{align}
Comparing coefficients of $\beta_1(z)^2$ between (\ref{phi-4-1}) and (\ref{phi-4-2}) yields that $c_2=0$.   Therefore $\beta_2=\beta_1^2$.   We shall henceforth drop the subscript and write $\beta$ in place of $\beta_1$ and $\beta^2$ in place of $\beta_2$.   

Having established that $\beta_2(z)=\beta(z)^2$ it is a straightforward matter to infer that $q_1=1$ by comparing coefficients between (\ref{F1}) and (\ref{F2}), as follows.  For each $n\geq1$, $\varphi_{2n}$ is the coefficient of $w^{2n}/(2n)!$ in (\ref{F1}), and the coefficient of $w^n/n!$ in (\ref{F2}).  Thus,
\begin{align}
\varphi_{2n}(z)&=(2n)!\left(c_{2n}+c_{2n-1}\beta(z)+\cdots+\frac{\beta(z)^{2n}}{(2n)!}\right)\label{2n-A}\\
&=n!\left(d_n+d_{n-1}\beta(z)^2+\cdots+\frac{\beta(z)^{2n}}{n!}\right)\label{2n-B}
\end{align}
and
\begin{align}
\varphi_{2n+2}(z)&=(2n+2)!\left(c_{2n+2}+c_{2n+1}\beta(z)+\cdots+\frac{\beta(z)^{2n+2}}{(2n+2)!}\right)\label{2n+2-A}\\
&=(n+1)!\left(d_{n+1}+d_{n}\beta(z)^2+\cdots+\frac{\beta(z)^{2n+2}}{(n+1)!}\right)\label{2n+2-B}
\end{align}
Comparing the coefficients of the first power of $\beta(z)$ in (\ref{2n-A}) with that in (\ref{2n-B}) reveals that $c_{2n-1}=0$.  Comparing constant terms between the same two formulas shows that 
\[
(2n)!c_{2n}=n!d_n
\]
while comparison of the coefficients of $\beta(z)^2$ between (\ref{2n+2-A}) and (\ref{2n+2-B}) reveals that 
\[
(2n+2)!c_{2n}/2=(n+1)!d_n.
\]
The unique solution to these latter two equations is $c_{2n}=d_n=0$.  This proves that $c_n=d_n=0$ for every $n\geq1$, whereby $q_1=1$ as claimed.   Thus 
\[
F_1(z,w)=e^{\alpha(z)+w\beta(z)}
\]
by Proposition~\ref{prop-10mod}. 
\end{pf}

Theorem~\ref{thm-Structure} may be proved using Proposition~\ref{prop-q1} as follows.  Since $F_1(z,w)=e^{\alpha(z)+w\beta(z)}$, comparison to (\ref{char}) shows that the moments of $A$ have the form 
\begin{equation}\label{psi-form}
\psi_n(z)=e^{\alpha(z)}\beta(z)^n,
\end{equation}
for $z\in\Omega$.  In particular, for all $z\in\Omega$,  $\psi_0(z)=e^{\alpha(z)}$ and $\psi_1(z)=e^{\alpha(z)}\beta(z)$.   Setting 
\begin{equation}\label{phi-meromorphic}
\varphi=\frac{\psi_1}{\psi_0},
\end{equation}
it follows that $\varphi$ is a meromorphic extension of $\beta$ from $\Omega$ to the disk $\ddd$.  Moreover, for every $z\in\Omega$,
\begin{equation}\label{extension}
\psi_n(z)=\psi_0(z)\varphi(z)^n
\end{equation} 
by equation (\ref{psi-form}).  Since $\Omega$ contains a condensation point, and both $\psi_n$ and $\psi_0\varphi^n$ are meromorphic in $\ddd$ for every $n\geq0$, the equation (\ref{extension}) extends to $\ddd$.  Furthermore, since $\psi_n$ is analytic, $\varphi$ cannot have any poles in $\ddd$, since otherwise $\psi_n$ too would have poles for sufficiently large $n$ by (\ref{extension}).  Thus $\varphi\in\sao$.   Since $\beta$ is non-constant (for $A$ of rank at least 2, by Proposition~\ref{prop-10mod}), (\ref{phi-meromorphic}) and part 1~of Proposition~\ref{prop-non-zero-zeroeth-moment} imply that $\varphi(\Omega)\subset\ddd$.   

Since $A$ is determined by its moments, it follows from (\ref{extension}) and the subsequent remarks that if $A$ has rank at least two, $A=M_\psi C_\varphi$, where $\psi=\psi_0\in\sao$, $\varphi\in\mathfrak{A}(\ddd)$ is non-constant, and $\varphi(\Omega)\subset\ddd$. The only remaining possibility for a non-zero operator $A\in\ssao$ is that it have rank 1, in which case it has the structure indicated in item 1~of Theorem~\ref{thm-Structure}.   Conversely, any operator corresponding to items 1~or 2~of Theorem~\ref{thm-Structure} is easily verified to belong to $\ssao$.   This completes the proof of Theorem~\ref{thm-Structure}.

\section{Conclusions\label{sec-conclusions}}

Sections \ref{sec-general} and \ref{sec-outer} illustrate the depth and scope of our main structure theorem as a bridge between classical P\'olya-Schur problems and analogous analytic problems.  On one hand, Theorem~\ref{thm-Structure} facilitates the resolution of a basic question in geophysical signal processing motivated by practical applications.  And on the other hand, the same theorem sheds new light on the basic theory of P\'olya-Schur problems in the classical setting.  

The results in \S\ref{sec-outer} make it possible to characterize explicitly the operators on $L^2(\real_+)$ that preserve the class of delayed minimum phase signals.  A minimum phase signal $f\in L^2(\real_+)$ is one that maximizes partial energy $\int_{0}^T|f(t)|^2\,dt$ among all functions having the same power spectrum as $f$, for all $T>0$.   Full details of this characterization and its relevance to seismic signal processing are laid out in \cite{GiLa:IP2012}.   

Several remarks are in order concerning the implications of Theorem~\ref{thm-General} for general P\'olya-Schur problems.  Firstly, the conclusion of the theorem is constructive, being formulated explicitly in terms of product-composition operators of the form $M_\psi C_\varphi$.  By contrast, consider \cite[Cor.~3]{BoBr:2009An}, which characterizes $\ssao$ in the case $\Omega=\ddd$ in algebraic terms as follows.  
\begin{thm}[From Corollary~3 in \cite{BoBr:2009An}]\label{thm-Cor}  A linear map $A:\complex[z]\rightarrow\complex[z]$ satisfies $A(\sss(\ddd))\subset\sss(\ddd)$ if and only if either:
\begin{enumerate}
\item there exist a linear functional 
$\nu:\complex[z]\rightarrow\complex$ and a polynomial $\psi\in\sss(\ddd)$ such that $A(f)=\nu(f)\psi$ for every $f\in\complex[z]$; or
\item letting $f_{w,n}(z)=(1+wz)^n$, the polynomial $Af_{w,n}$ is $\ddd$-stable for every $w\in\ddd$ and $n\geq0$.   
\end{enumerate}
\end{thm}
In other words, letting $\mathcal{M}$ denote the set of polynomials
\[
\mathcal{M}=\left\{ f_{w,n}\,\left|\, w\in\ddd\mbox{ and }n\geq0\right.\right\},
\]
Theorem~\ref{thm-Cor} reduces the case $\Omega_1=\Omega_2=\ddd$ to the non-trivial problem of determining all linear operators $A:\complex[z]\rightarrow\complex[z]$ such that 
$
A(\mathcal{M})\subset\sss(\ddd).
$
Theorem~\ref{thm-General} solves the latter problem explicitly, thus providing complementary information in cases where there is a known characterization.  

Apart from giving constructive solutions in cases where there are known, non-constructive characterizations, Theorem~\ref{thm-General} solves the previously open case of $\Omega_1=\Omega_2=\overline{\ddd}$, as well as the array of cases where $\Omega_1$ and $\Omega_2$ do not correspond to circular regions or their boundaries.  

A final remark concerns the general P\'olya-Schur problem of characterizing $\ssao$ for $\Omega\subset\complex$.  Theorem~\ref{thm-General} implies that if $\Omega$ is bounded with non-empty interior, then $\ssao$ consists only of rank 1 operators and product-composition operators.  It is clear, for example, that if $\complex\setminus\Omega$ is convex then $\ssao$ includes differential operators (by the Gauss-Lucas Theorem), which are not rank 1 and are not expressible as product-composition operators.   But there is a substantial gap between non-convexity of $\complex\setminus\Omega$ and $\Omega$ being bounded with non-empty interior, and it remains an open question to characterize precisely all regions $\Omega\subset\complex$ such that $\ssao$ consists exclusively of rank 1 operators and product-composition operators.   

\vspace*{3ex}
\noindent\emph{Acknowledgements.}
Research for the present paper was supported by NSERC Discovery Grants and the MITACS project POTSI.  The first author wishes to the thank the Mathematics Department at PUC-Rio for providing warm hospitality and a stimulating research atmosphere during his visit there in 2011.  Thanks also to the anonymous referees for their comments.

\end{document}